\documentclass{birkjour}
\usepackage[latin1]{inputenc} 
\usepackage{graphicx}
\usepackage{amssymb}
\usepackage{fontawesome}
\usepackage{amsthm}
\usepackage{amsmath}
\usepackage{url}
\usepackage{tikz}
\usetikzlibrary{decorations.markings}
\tikzstyle{point}=[circle, draw, inner sep=0pt, minimum size=6pt]
\usepackage{caption}
\theoremstyle{plain} 
\newtheorem{thm}{Theorem}[section]

\theoremstyle{definition} 
\theoremstyle{remark} 

\theoremstyle{definition} 

\newtheorem{rem}[thm]{Remark}

\begin{document}
	\title[Area-filling curves]
	{Area-filling curves} 
	
	\author[Maria Chiara Nasso]{Maria Chiara Nasso}
	
	\address{Dipartimento di Ingegneria Informatica, Modellistica, Elettronica e Sistemistica, Universit\`a della
		Calabria,\\ Rende, CS, Italy.}
	
	\email{mc.nasso@dimes.unical.it}

	\author{Aljo\v{s}a Vol\v{c}i\v{c}}
	\address{Dipartimento di Matematica e Informatica,\\ Universit\`a della Calabria,\\ Rende, CS, Italy.}
	\email{aljosa.volcic@unical.it}
	

	\maketitle

\begin{abstract} 
	In this paper we study area-filling curves, i.e. continuous and injective mappings defined on $[0,1]$ whose graph has positive measure. Current literature calls them ``Osgood curves", but their invention is due to H. Lebesgue. Stromberg and Tseng constructed {\it homogeneous} area-filling curves and offered an elegant example. We show that an appropriate variant of Knopp's construction attains the same homogeneity result. In Section 4 we discuss briefly the existence of an ``invasive" curve, i.e. a continuous and injective mapping from the half-open interval $[0,1[$ to the unit square, whose image has measure 1. In the last section we discuss several aspects of the Lance-Thomas curve, connecting it with the other  construction due to Stromberg and Tseng.\\
\end{abstract}

\section{Introduction}
This paper is dedicated to some questions raised by Sagan's entertaining book ``Space-filling curves" \cite{Sagan1}
and the deep and analytically highly non trivial paper \cite{Stromberg}.
\\
Peano, Hilbert (\cite{Peano,Hilbert}) and other authors constructed examples of 
{\it surjective} continuous mappings $f:[0,1]\rightarrow [0,1]^2$ which are called {\it space-filling curves}. It is well-known that such curves cannot be injective by a result of Netto \cite{netto}. \\
Nowadays, space-filling curves have several useful and surprising applications in global optimization (see e.g. \cite{Sergeyev}).\\
In analogy we shall call {\it area-filling curves} those continuous and {\it injective} mappings from $[0,1]$ to $[0,1]^2$ which cover a set of positive measure. \\
To our knowledge, there are ten constructions of area-filling curves. Moreover there are two existence theorems. 
Le\-bes\-gue and Osgood proposed two constructions which are in essence equivalent. The others are due to Sierpi\'nski \cite{Sierp},
Knopp \cite{Knopp}, 
Gelbaum and Olmsted (who offered two constructions) 
\cite{GelOlm}, 
Lance and Thomas \cite{Lance}
(with the variant due to Sagan 
\cite{Sagan2}
 and another one due to Ravishankar and Salas \cite{Ravishankar}), and two different constructions are due to Stromberg and Theng \cite{Stromberg}. The last paper contains also a non-constructive existence proof suggested to the authors by Sadahiro Saeki. Another interesting non-constructive existence proof was given in \cite{BalKh} by Balcerzak and Kharazishvili using the Denjoy-Riesz theorem. All the results we discuss are variant of the following\\

	\begin{thm}
		Given a regular set $Q$ having measure $1$,  for every $\beta \in ]0,1[$ there exists a continuous injective curve $\gamma_\beta:[0,1]\rightarrow Q$ such that $\lambda_2(\gamma_\beta([0,1]))$ $=\beta.$
	\end{thm}
	
		We indicate with $\lambda_m$ the m-dimensional Lebesgue measure.
		In \cite{Sierp}
		and \cite{Knopp} 
		  $Q$ is a triangle, in \cite{Leb2}
		   it is a sector of a circular crown. In all the other constructions $Q$ is the unit square.
	Sierpi\'nski's construction has the drawback that it does not allow to chose {\it an arbitrary} $\beta\in ]0,1[$.\\
	Current literature consistently calls such curves ``Osgood curves", but we suggest to call them ``area-filling curves", since we will show in Section 2 that Lebesgue invented such curves one year before Osgood. On the other hand we cannot call them Lebesgue curves, since there exists already a space-filling curve named after him.
	 
	A curve is a continuous mapping from an interval $I\subset \mathbb{R}$ into $\mathbb{R}^m$. In our paper $I$ will be a bounded closed interval, except in Section 4, where topological circumstances will suggest us to take it half-open.
	
	If $\varphi$ and $\psi$ are two curves defined on $[a,b]$ and $[a',b']$, respectively, with values in $\mathbb{R}^m$, we will call them {\it equivalent} if there exists an increasing homeomorphism $s(\cdot)$ from $[a,b]$ to $[a',b']$ such that
	$\varphi(\cdot)=\psi(s(\cdot))$ for every $t\in [a,b]$.

Sometimes it is convenient to identify curves with their equivalence classes, while in other cases it is important to distinguish between two equivalent curves.

This ambivalent approach can be seen already in the papers which are at the origin of space- and area- filling curves: Peano provides the analytic expression of his curve (and shows no picture), while Hilbert only shows pictures, saying that they lead trivially to their parametric representation.

Area-filling curves are not part of the traditional program offered to students of mathematics in spite of the simplicity of some constructions and the
influential mathematicians who contributed to their study.  Area-filling curves tend to be forgotten and rediscovered (see 
\cite{Sagan1}, pp. 132-133). 

We think, on
the contrary, that they should belong to the permanent display of amusing objects aimed to arouse the interest of students (like the Peano and Hilbert curve, the snowflake curve,
the Cantor set and the Cantor function, the Sierpi\'nski triangle, the M\"{o}bius strip, Klein's bottle and the rest of the zoo).
The origin itself of the area-filling curves is curious. In Section 2 we will discuss their early history.
In Section 3 we will discuss the {\it homogeneous constructions} of area-filling curves, a concept introduced in \cite{Stromberg}. We will show that an appropriate (and simpler) variant of Knopp's construction also has the homogeneous property.
In Section 4 we will discuss a related question, the existence of a continuous and injective mapping 
$\varphi:[0,1[\rightarrow \mathbb{R}^m$ such that $\lambda_m(\varphi([0,1[))=1$. 
In the last section we will generalize the construction of Lance and Thomas \cite{Lance} and construct a class of area-filling curves which depend (for any given $\beta$) on continuously many parameters. They will be parametrized on symmetric Cantor sets which will turn out to be homeomorphic to the {\it essential trace} of the curve.

\section{Early history}

H. Lebesgue, in his famous thesis, considered a measure-theoretic problem which can be reformulated in the following way:
{\it Suppose $A$ is the interior of a Jordan curve. Is $A$ Peano-Jordan-measurable?}\\
Put it in other words, the question is if there are Jordan curves whose image (sometimes called trace) has positive Jordan measure.

This was unknown at that time and Lebesgue sketched in a footnote of his thesis (page 17), the general idea for the proof of their existence. It was based on the construction of the Peano curve.

Osgood \cite{Osgood}, one year after the publication of Lebesgue's thesis, constructed independently an example of an injective curve having positive area, based on the construction of the Peano curve.

In the same year \cite{Leb2} Lebesgue returned to his measure-theoretic question, elaborating a particular instance of his general idea and proposed a similar construction, using again Peano's idea.

It is curious that nobody noticed in this 120 years that Lebesgue in \cite{Leb2} (quoted in  \cite{Stromberg}, for instance) refers to \cite{Leb}. Also the little measure-theoretic question, which motivated it, seems to have been forgotten.

\section{Homogeneity}

Sierpi\'nski and Knopp criticized the fact that in Osgood's construction some arcs have area zero and they proposed two similar constructions which eliminate this drawback, constructing two different examples of
 area-filling curves $\gamma$ which have only arcs of strictly positive area. We will call curves with this property {\it positive} curves. Another positive curve has been constructed by Gelbaum and Olmsted (their first example).

Stromberg and Tseng \cite{Stromberg} presented a construction of a positive area-filling curve which has an additional property.

They wrote: ``there is a `homogeneity' feature which is desirable but not guaranteed by Knopp's and Siepi\'nski's construction: given $\beta \in ]0,1[$, one may ask that for every subset E of $[0,1]$ the image $\gamma(E)$ have area equal to $\beta$ times the `length' of $E$.\\

Here it is important to distinguish between curves as functions 
 and curves as equivalence classes as we will see.
 
But first we want to show that an appropriate variant of the Knopp construction, which seems to us simpler, provides also this homogeneity feature.

\begin{thm}
	For every $\beta \in ]0,1[$ there exists a homogeneous construction of a Knopp curve with $\lambda_2([0,1])=\beta$.
\end{thm}
\begin{proof}
	Take a triangle $\triangle ABC$ having area 1. At the first step subtract from it a triangle with vertex in $C$ and base $DE$ on the segment $AB$ to obtain two closed triangles $T_0=\triangle ADC$ and $T_1=\triangle EBC$  having equal area (details will be given later). See Figure \ref{Knoppdis}.
	
	\begin{figure}[!h]
		
		\centering
		
		\includegraphics[width=0.7\textwidth]{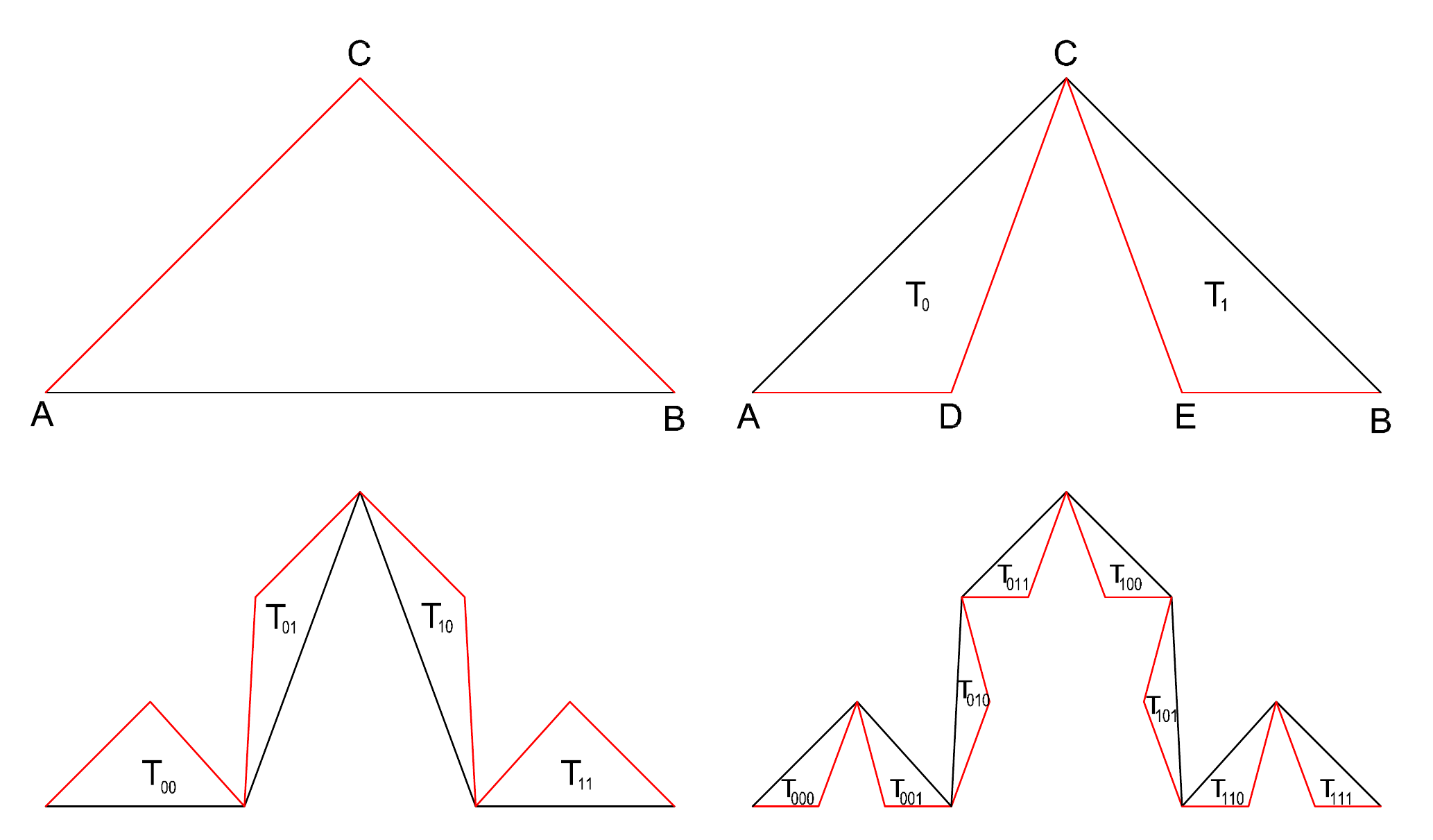}
		
		\caption{First four steps of the generation of Knopp homogeneous area-filling curve.}\label{Knoppdis}
	\end{figure}
	
	Repeat now this construction subtracting from $T_0$ a triangle with vertex in D and base a segment on $AC$, and from $T_1$ a triangle with vertex in E and base a segment on $BC$. This can be done obtaining a chain of four triangles  (so that they have equal area) $T_{00}$, $T_{01}$, $T_{10}$ and $T_{11}$, which connects $A$ to $B$.
	
	Iterating this procedure, we get nested chains $V_n$ of triangles having the same area connecting $A$ to $B$. The indices of $T_{\varepsilon_1, \varepsilon_2, \dots, \varepsilon_n}$, with $\varepsilon_k \in \{0,1\}$, are chosen so to establish an ordered correspondence between the triangles and the binary $n$-digit fractions
	$0.\varepsilon_1 \varepsilon_2, \dots, \varepsilon_n$.
	
	$V=\bigcap_{n=1}^{\infty} V_n$ is a compact connected set and for any $\beta \in  ]0,1[$, it is possible to chose (in many ways) the areas of the subtracted triangles so that $\lambda_2(V)=\beta$.

	The crucial step in Knopp's construction is to prove that the diameters of the triangles $T_{\varepsilon_1, \varepsilon_2, \dots, \varepsilon_n}$ tend to zero. This technical step proves that there exists a continuous and injective mapping $\gamma$ from $[0,1]$ onto $V$ which assigns to every real number written in the binary form $0.\varepsilon_1 \varepsilon_2, \dots, \varepsilon_n, \dots$ the uniquely determined point $v\in V$ which belongs to the intersection $\bigcap_{n=1}^{\infty}  T_{\varepsilon_1, \varepsilon_2, \dots, \varepsilon_n}$. 
	
	Our choice of the triangles is different, so we have to prove that our construction produces triangles whose diameters tend to zero, so we can use the last part of Knopp's proof.
	
	Let $\beta \in ]0,1[$ and take a sequence of real numbers $\{r_k\}$
	such that  $r_k\in ]0,1[$ and $ \prod_{k=1}^{\infty}(1-r_k)=\beta$.\\
	At the first step let $T_0$ and $T_1$ have area $\frac{1-r_1}{2}$. \\
	At the next step let $T_{00}$, $T_{01}$, $T_{10}$ and $T_{11}$
	have area $\frac{(1-r_1)(1-r_2)}{4}$.\\
	In general, let the $2^n$ triangles of $V_n$ have area $\frac{1}{2^n}\prod_{k=1}^n(1-r_k)$.

	Since the sequence $\{r_k\}$ tends to zero, there exists $n_0$ such that $r_n<\frac{1}{8}$ for $n\ge n_0$. Take an $n\ge n_0$ and consider a triangle in the chain $V_{ n}$  with vertices N, P and O. 
	
	We are going to show that in three steps the diameters of the triangles are reduced by a factor at most $\frac{3}{4}$.\\
	The triangle $\triangle NQR$ (see Figure \ref{ksegnato23} a)) 
	is similar to $\triangle NPO$ with similarity ratio $\frac{1}{2}$. Since $\triangle NQ_1R_1\subset \triangle NQR$, this triangle and all the triangles of the successive chains will have diameter less than  $\frac{1}{2}$ times the diameter of $\triangle NPO$. The same argument applies to $\triangle Q_2PS_1$.
	
	\begin{figure}[!h]
	
	\centering
	
	\includegraphics[width=0.90\textwidth]{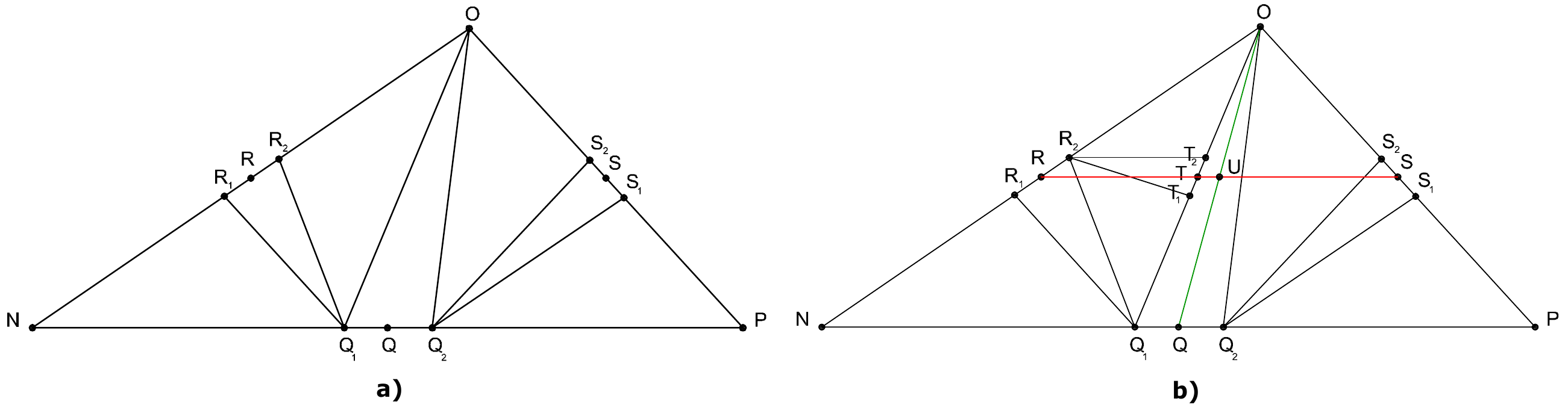}
	
	\caption{Triangle of $V_n$ and construction at $(n+2)$th step (left) and $(n+3)$th step (right).}\label{ksegnato23}
\end{figure}
	Let us now turn to $\triangle R_2Q_1O$ and the operations of the next step (see Figure \ref{ksegnato23} b)).
$\triangle R_2Q_1O$ generates $\triangle R_2T_2O$ and $\triangle R_2T_1Q_1$. The first one is contained in $\triangle ORS$ and its diameter is smaller than one half of the diameter on $\triangle NPO$.
	$\triangle R_2Q_1T_1$ offers greater resistance. But it is easy to see that it is contained in the image of 
	$\triangle NPO$ under a homothety which shrinks it by a factor $\frac{3}{4}$ and keeps $N$ fixed. This is due to the fact that $r_n < \frac{1}{8}$. 
	
	Now Knopp's argument applies and it follows that $V$ is the image of a continuous injective mapping.\\
	Let us prove now the homogeneity property of $\gamma$.
	Let $T=T_{\varepsilon_1, \dots, \varepsilon_n}$ be a triangle of the chain $V_n$ and denote by $X$ and $Y$ its vertices which connect it to the remaining part of the chain. Possibly $X=A$ or $Y=B$. The points $A$, $X$, $Y$ and $B$ are in that order on $V$.
	We have that
	$\{X\}=\bigcap_{m=0}^{\infty}T_{\varepsilon_1, \dots, \varepsilon_n0\dots 0}\,,
	$
	where $\varepsilon_n$ is followed by $m$ $0$'s.
	Therefore $X=\gamma(0.\varepsilon_1 \dots \varepsilon_n)$.
	
	On the other hand
	$\{Y\}=\bigcap_{m=0}^{\infty}T_{\varepsilon_1, \dots, \varepsilon_n1\dots 1}\,,
	$
	where $\varepsilon_n$ is followed by $m$ $1$'s.
	Therefore $Y=\gamma(0.\varepsilon_1 \dots \varepsilon_n, \bar 1)$.
	
	It follows that the interval $I\subset [0,1]$, mapped by $\gamma$ onto the arc
	$\overset{\frown}{X Y}$, has length $\frac{1}{2^n}$.
	
	On the other hand $T$ has area 
	$$\prod_{i=1}^n\frac{(1-r_i)}{2}=\frac{1}{2^n}\prod_{i=1}^n(1-r_i)\,.$$
	The $2^m$ triangles belonging to $V_{n+m}$  and contained in $T$ have total area 
	$$2^m\frac{1}{2^{n+m}}\prod_{i=1}^{n+m}(1-r_i)=\frac{1}{2^{n}}\prod_{i=1}^{n+m}(1-r_i)\,.
	$$
	Taking the limit for $m\rightarrow \infty$, it follows that 
	$$\lambda_2(\overset{\frown}{X Y})=\frac{1}{2^{n}}\prod_{i=1}^{\infty}(1-r_i)=\beta \frac{1}{2^{n}}=\beta \lambda_1(I)\,.$$
	The equation
	$\lambda_2(\gamma(E))=\beta \lambda_1(E),$
	for every measurable set $E$ contained in $[0,1]$, follows now from the density of binary intervals and their finite unions among measurable sets.
\end{proof}

We may speak about a homogeneous parametrization of a curve. But curves alone cannot be considered homogeneous and we will show that any {\it positive} curve can be reparametrized so to get a homogeneous parametrization.

\begin{thm}
	Let $\gamma$ be a positive area-filling curve parametrized on $[0,1]$ (with $\lambda_2(\gamma([0,1]))=\beta$). Then there exists an increasing homeomorphism $h$ from $[0,1]$ to $[0,1]$ such that $\gamma_1= \gamma \circ h$ is homogeneously parametrized.
\end{thm}

\begin{proof} The function 
	$h^{-1}(t)=\beta^{-1}\lambda_2(\gamma([0,t]))
	$
	is an increasing homeomorphism from $[0,1]$ to $[0,1]$. Let us denote its inverse with $h$.
	
	Since $h^{-1}(h(s)) = s=\beta^{-1}\lambda_2(\gamma([0,h(s)])$,
	we have
	$$\lambda_2(\gamma_1(s))=\lambda_2(\gamma([0,h(s)])=\beta \,s\,.$$
	as required.
\end{proof}

On the other hand, if we take the homogeneous representation of the Knopp's curve just constructed and consider the increasing homeomorphism $s(\cdot): [0,1] \rightarrow [0,1]$ $s(t)=t^2$, we notice that $\gamma_1(t)=\gamma(t^2)$ is such that $\lambda_2(\gamma_1([0,\frac{1}{2}]))=\frac{1}{4}\beta$ and $\lambda_2(\gamma_1([\frac{1}{2},1]))=\frac{3}{4}\beta$, so $\gamma_1$ is not homogeneously parametrized.

\section{Invasive curves}
The natural question whether there exists an injective continuous curve $\gamma$ with image in $[0,1]^2$ such that $\lambda_2(\gamma([0,1]))=1$ has an immediate negative answer: $\gamma([0,1])$ is compact and the only compact subset of $[0,1]^2$ with measure 1 is the square itself. But this is excluded by Netto's result.\

We will mention here two results related to this question which can be deduced or found in \cite{Stromberg} and \cite{PietroePaolo}, respectively.

Sadahiro Saeki contributed to \cite{Stromberg} with a theorem which is not constructive, since it uses the deep result due to von Neumann, Oxtoby and Ulam, the Homeomorphic Measures Theorem (HMT) (see \cite{oxeul} and \cite{GofPed}). He proved just the existence of non-intersecting area-filling curves.
Looking at his proof, however, it is easy to draw the following (non-constructive) conclusion.

\begin{thm}
There exist sequences $\{\gamma_n\}$ of non-intersecting injective curves defined on $[0,1]$ such that
$\sum_{n=1}^{\infty}\lambda_2(\gamma_n([0,1]))=1\,.$
\end{thm}

A result which gets closer to an answer to the question posed at the beginning of the section is contained in \cite{PietroePaolo}. The Authors do not state it explicitly, but it can be found in the discursive part of their article (p. 176, lines -6, -5):
\begin{thm}
There exists a continuous and injective mapping $\gamma$ from $[0,1[$ to $[0,1]^2$ such that $\lambda_2(\gamma([0,1[))=1$.	
\end{thm}
It is an easy consequence of their nice main result, the Thread Theorem.

\section{Lance-Thomas curves}
In the introduction we mentioned non constructive existence proofs. One is contained in \cite{BalKh}. It is not explicitly stated, but is described in the text of the proof of their Theorem 2.2. It uses the Denjoy-Riesz theorem, which asserts that

\begin{thm}
	Any bounded totally disconnected set in the plane can be covered by a Jordan arc.
\end{thm}

The existence of an area-filling curve follows easily from the previous result: fix $\beta>0$, take a planar totally disconnected set of measure $\beta$ contained in $[0,1]^2$ and cover it with a Jordan arc. So one obtains an area-filling curve of measure at least $\beta$ (since the Jordan arc may contain other area-filling arcs). 
The easiest way of producing a totally disconnected compact set in the plane of positive measure is to take a Cantor set $P\subset [0,1]$ of positive measure and to consider then $P\times P$. 
Having the Denjoy-Riesz theorem in mind, the proofs of Lebesgue, Osgood and the second construction by Gelbaum and Olmsted can be immediately related to this scene, but there is no doubt that the direct inspiration came from Peano's and Hilbert's set. 

It should be noted however that the first of the two constructions from \cite{Stromberg} can also be related to this trace: fixed a $\beta >0$, they let  $\alpha=\sqrt{\beta}$ and constructed a symmetric Cantor set $P\subset [0,1]$ having measure $\alpha$. Then they exhibited explicitly segments which connect points of $P\times P$ and form a Jordan arc covering $P\times P$. Its measure is $\beta$, since the segments have measure zero.
The construction is explicit, but very complicated (it takes pages 34-43) and it does not present any evident symmetry. They do not provide any picture of their set and also our attempts to obtain a suggestive figure failed. 

 It should be noted that \cite{BalKh} has been published five years later, so there is no direct connection between the two papers, and we just imagine how the authors have argued.

We will now make a change to the construction of the Lance-Thomas curve to fit into this line of thought.
Let $\beta \in  ]0,1[$ and consider a sequence $\{a_n\}$ of real numbers,
called the reduction coefficients, such that $0 < a_n < 1$ and
$\lim_{n \rightarrow \infty}(a_1 \cdot a_2 \cdot \dots \cdot a_n)^2=\beta \,.$

Our construction of the Lance-Thomas curve differs from the previous ones (\cite{Lance}, \cite{Sagan1} and \cite{Sagan2}) by the parametrization. Note that Sagan parametrized it on the ternary Cantor set and did not follow \cite{Lance} who simply divided $[0,1]$ into seven equal parts.

\begin{figure}[!h]
	
	\centering
	
	\includegraphics[width=0.6\textwidth]{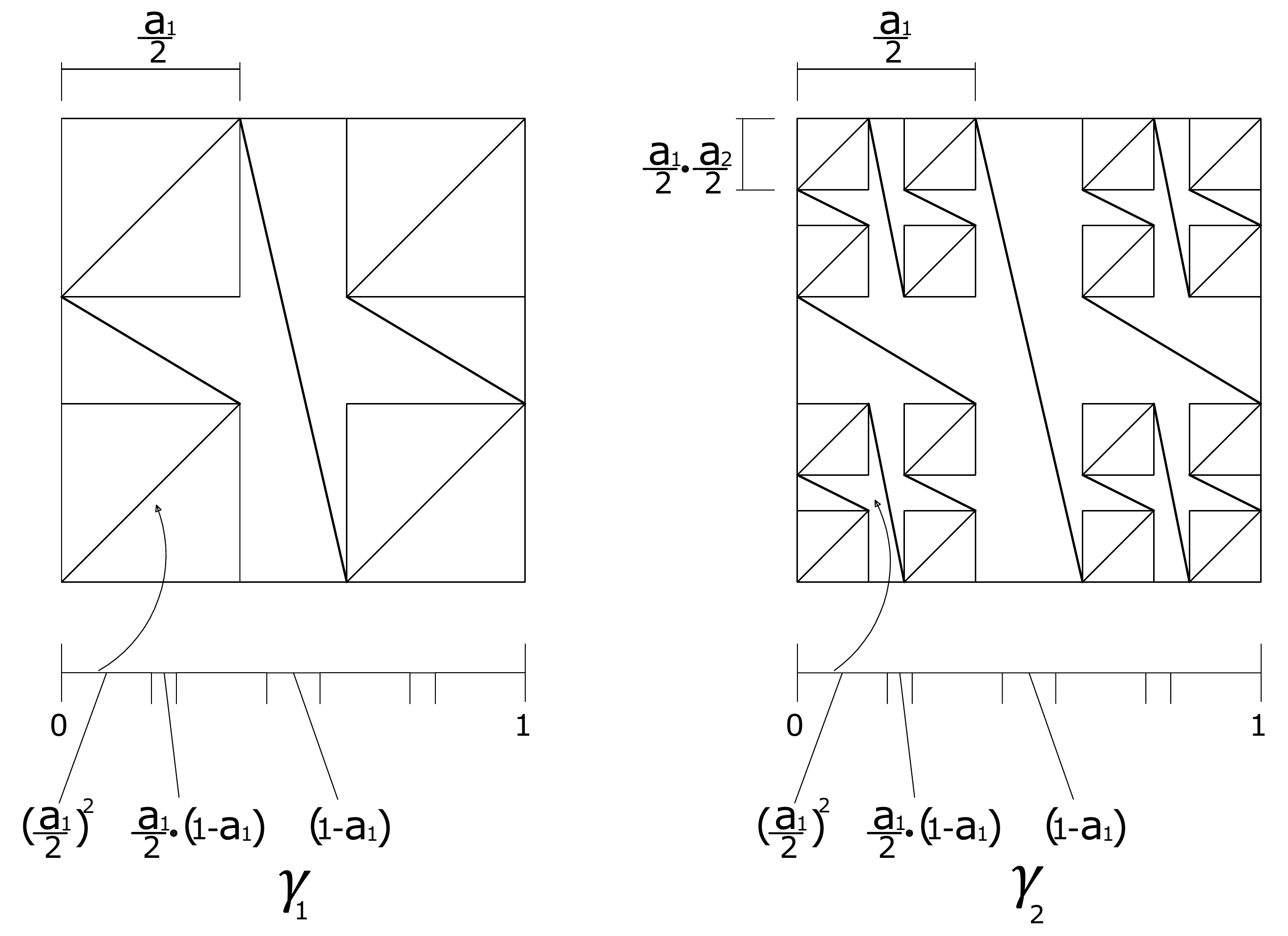}
	
	\caption{{Construction of $\gamma_1$ and $\gamma_2$.}}\label{fnNostra}
\end{figure}

To construct the first approximation, consider (following the original idea) four disjoint squares, having side length $\frac{a_1}{2}$ , placed at the corners of $[0, 1]^2$ and joined by three segments (joints) as in Figure (\ref{fnNostra}). 
The curve $\gamma_1$ is obtained adding to the three joints the diagonals of the squares to get a polygonal line which connects the left lower corner with the right upper corner of the square. 
Denote by $A_1$ the union of those four closed squares.

We choose to subdivide $[0, 1]$ into seven subintervals so that the first, the third, the fifth and the seventh have length $(\frac{a_1}{2})^2$. Denote by $C_2$ the union of these four intervals. The second and sixth intervals have length $\frac{a_1}{2}(1-a_1)$, while the central has length $(1-a_1)$. The function $\gamma _1$ maps linearly the even intervals (second, fourth, sixth) in the segments that connect the squares of $ A_1$ and the remaining subintervals in the natural order into the diagonals of the squares as shown in Figure (\ref{fnNostra}).

There is no mystery in the choice of the length of the intervals: they are equal to the areas of the seven rectangles crossed by the seven parametrized segments.

The curve $\gamma_2$ is constructed analogously, operating now in the four squares of the first generation using the reduction coefficient $a_2$.
Denote by $A_2$ the set made up of sixteen squares obtained
putting into each of the four squares of $A_1$ four smaller squares, placed at the vertices, of side length $\frac {a_1}{2}\frac{a_2}{2}$.  Trace their junctions as in Figure (\ref{fnNostra}). 

The four intervals of $C_2$ are subdivided as in the first step using the reduction coefficient $a_2$.
Each interval of $C_2$ is subdivided into seven parts: the first, third, fifth and seventh have length $(\frac {a_1}{2}\frac{a_2}{2})^2$. 

Let us denote by $C_4$ the union of these four disjoint closed subintervals. The second and sixth have length $(\frac{a_2}{2})^2 \cdot\frac{a_2}{2}(1-a_2) $ and the fourth has measure $(\frac{a_2}{2})^2 \cdot (1-a_2) $.

The function $\gamma_2$ maps linearly the intervals of $C_4$ into the diagonals of the $A_2$ squares and the remaining subintervals are mapped in the junction segments, as in Figure (\ref{fnNostra}), forming so a connected polygonal.

We iterate the construction, putting into each of the squares of $A_{n-1}$ four squares and three junctions as in step one and two, using at each step the reduction factor $a_n$. Iterating, we have $\lambda_2(A_n)=(a_1\cdot a_2 \dots a_n)^2$ and the set $C_{2n}$ (which is the $2n^{th}$ step of the construction of a symmetric Cantor set $C$) consists of $4^n$ disjoint closed intervals having one-dimensional measure equal to the two-dimensional measure of the rectangle having as diagonal the segment in which the interval is mapped, so $\lambda_1(C_{2n})=(a_1\cdot a_2 \dots a_n)^2$.

\noindent

Note that $\gamma_{n+1}=\gamma_n$ on the complement of $C_{2n}$. On the other hand, for every $k$ and each $x \in C_{2n}$, $||\gamma_{n+k}(x)-\gamma_n(x)||\leq \left( \frac{a_1}{2}\cdot\frac{a_2}{2}\cdot...\cdot\frac{a_n}{2}\ \right)\cdot \sqrt{2}$, which tends to $0$. Thus the sequence $\{\gamma_n\}$ converges uniformly to a continuous curve $\gamma$ such that 
$$\lambda_2(\gamma([0,1]))=\lambda_2\left(\bigcap_nA_n\right)=\lim_{n \rightarrow \infty}(a_1\cdot a_2 \dots a_n)^2=\beta\,.$$

\begin{rem}\label{misurapreservata}
	If we denote with $B_n$ the set which consists of the union of squares of $A_n$ for each n, the restriction of the limiting curve $\gamma$ to $C$ is a homeomorphism between $C$ and $B=\cap_n B_n$ that preserves the measure. We call the compact set $B$ the {\it essential image} of $[0, 1]$ under $\gamma$.

	It is easy to see that for any Borel set $E\subset [0, 1]$, $\lambda_2(\gamma(E\cap C))=\lambda_1(E\cap C)$, so our parametrization is sort of a ``homogeneous" one.
	
	Note also that $B$ is a cartesian product of two symmetric Cantor sets contained in $[0,1]$ and that the curve $\gamma$ is a Jordan arc containing $B$. So this construction is an illustration of the Riesz-Denjoy theorem.
\end{rem}


\begin{thebibliography}{}
	\bibitem{BalKh} M. Balcerzak, A. Kharazishvili, \textit{On Uncountable Unions and Intersections of Measurable Sets}, Georgian Math. J. 6 (1999), no. 3, 201--212.
	\bibitem{Falconer} K. J. Falconer, \textit{The Geometry of Fractal Sets}, Cambridge University Press, 1985.
\bibitem{Frechet} M. Fr\'echet, \textit{Sur quelques points du calcul fonctionnel},  Rend. Palermo { 27}, (1906), 1--74.	
\bibitem{GelOlm} B. R. Gelbaum, J. M. H. Olmsted, \textit{Counterexamples in Analysis}, Holden-Day, 1964.
\bibitem{GofPed} C. Goffman, G. Pedrick, \textit{A proof of the homeomorphism of Lebesgue-Stieltjes measure with Lebesgue measure}, Proc. Amer. Math. Soc. 52 (1975), 196--198.
\bibitem{Hilbert} D. Hilbert, \textit{\"{U}ber die stetige Abbildung einer Linie auf ein Fl\"{a}chenst\"{u}ck}, Math. Ann. 38 (1891), no. 3, 459--460.
\bibitem{Knopp} K. Knopp, \textit{Einheitliche Erzeugung und Darstellung der Kurven von Peano}, \textit{Osgood und von Koch}, Arch. Math. Phys. no. 26, (1917), 103-115.
\bibitem{PietroePaolo} P. Haj{\l}asz, P. Strzelecki, How to measure volume with a thread, Amer. Math. Monthly {112}, (2005), 176-179.
\bibitem{Lance} T. Lance and E. Thomas, \textit{Arcs with Positive Measure and a Space-Filling Curve},  Amer. Math. Monthly 98 (1991), no. 2, 124--127.
\bibitem{Leb} H. Lebesgue, \textit{Int\'egrale, longueur, aire}, These, 30 juin 1902, Facult\'e des Sciences de Paris.
\bibitem{Leb2} H. Lebesgue, \textit{Sur le probl\'eme des aires}, Bull. Soc. Math. France 31 (1903), 197--203.
\bibitem{netto} E. Netto, \textit{Beitrag zur Mannigfaltigkeitslehre}, J. Reine Angew. Math. 86 (1879), 263--268.
\bibitem{Osgood} W. F. Osgood,  \textit{A Jordan Curve of Positive Area}, Trans. Amer. Math. Soc. 4 (1903), no. 1, 107--112.
\bibitem{oxeul} J. C. Oxtoby and S. M. Ulam, \textit{Measure-Preserving Homeomorphisms and Metrical Transitivity}, Ann. of Math. (2) 42 (1941), 874--920.
\bibitem{Peano} G. Peano, \textit{Sur une courbe, qui remplit toute une aire plane}, Math. Ann. 36 (1890), no. 1, 157--160.
\bibitem{Ravishankar}K. Ravishankar and H. Salas,\textit{ On the existence of locally heavy arcs}, Rev. Un. Mat. Argentina 36 (1990), 101--110 (1992). 
\bibitem{Sagan2} H. Sagan, \textit{A Geometrization of Lebesgue's Space-Filling Curve}, Math. Intelligencer 15 (1993), no. 4, 37--43.
\bibitem{Sagan1} H. Sagan, Space-filling curves, Universitext, New York, 1994.
\bibitem{Sergeyev} Y. D. Sergeyev, R. G. Strongin and D. Lera, \textit{Introduction to Global Optimization Exploiting Space-Filling Curves}, Springer, 2013.
\bibitem{Sierp} W. Sierpi\'{n}ski, \textit{Sur une courbe non quarrable}, Bull. Acad. Sci. de Cracovie {Sci. math. et nat., S\'{e}rie A} ({1913}), 254-263. 
\bibitem{Stromberg} K. Stromberg and S. Tseng, \textit{Simple plane Arcs of positive Area}, Expo. Math. 12 (1994), 31-52.
 \bibitem{Vallin} R. W. Vallin, \textit{The Elements of Cantor Sets with applications}, 2013.
\end{thebibliography}
\end{document}